\documentclass[a4paper]{article}
  \usepackage[dvips]{graphicx}
\usepackage[cmex10]{amsmath}
\usepackage{amsfonts,amssymb,amsthm}
\usepackage[font=footnotesize]{subfig}
\newtheorem{thm}{Theorem}
\newtheorem{lem}{Lemma}
\newtheorem{defn}{Definition}
\newcommand{\myp}[1]{\mathbf{\Phi}[#1]}

\newcommand{\myce}[2]{\mathbb{E}\{#1|#2\}}

\begin{document}
%
\title{Scheduling Under Fading and Partial Channel Information}

\author{
Santanu Mondal and Vinod Sharma\\
Department of Electrical Communication Engineering,\\
Indian Institute of Science, Bangalore -- 560 012, India.\\
Email: \{santanu,vinod\}@ece.iisc.ernet.in
}


%


This work has been submitted to the IEEE for possible publication. Copyright may be transferred without notice, after which this version may no longer be accessible.

\date{}
\maketitle

\begin{abstract}
We consider a scheduler for the downlink of a wireless channel when only partial channel-state information is available at the scheduler. We characterize the network stability region and provide two throughput-optimal scheduling policies. We also derive a \emph{deterministic} bound on the mean packet delay in the network. Finally, we provide a throughput-optimal policy for the network under QoS constraints when real-time and rate-guaranteed data traffic may be present.
\end{abstract}


%

\section{Introduction}

 

Scheduling has always been an indispensable part of resource allocation in wireless networks. 
The seminal work of Tassiulas \emph{et al.} \cite{tas1}, and later \cite{tas2}, \cite{tas3} considered the case where both channel states and queue lengths are fully available to the scheduler. It was shown that the MaxWeight algorithm, which serves the longest connected queue, is throughput-optimal.
Subsequently, the MaxWeight algorithm was found to be throughput-optimal in many other settings as well (\cite{sto2}-\cite{shak2} and the references therein) using tools from Lyapunov optimization. Some other works (see \cite{R1}, \cite{R2}, \cite{ness2}) also approach the scheduling problem using convex optimization and dual decomposition techniques. \cite{sri3} even considers the role of imperfect queue length information on network throughput, showing that the stability region does not reduce. But, in all these cases, accurate information about channel-state is assumed as a modeling simplification.

In a real-life network, \emph{e.g.}, Long Term Evolution (LTE) \cite{lte} or IEEE 802.16e WiMAX, the channel-state information fed back to the transmitter can have uncertainty. The primary reason being that although resource-allocation is done at the \emph{finer granularity} of a Physical Resource Block (PRB), channel-state information is still fed back at the \emph{coarser granularity} of a subband, which is a group of PRBs. This is done to reduce the feedback traffic from the users to the Base Station (BS). However, this \emph{averaging} causes information loss and hence, the resulting uncertainty at the scheduler. Moreover, uncertainty might be present in the channel-estimates because of the very process of estimation. 

Some recent works have, hence, tried to model this uncertainty in the channel-estimate. In \cite{kar1}, the authors show that infrequent channel-state measurement, unlike infrequent queue length measurement, reduces the maximum attainable throughput. \cite{neelesh} considers the effect of inaccuracy of channel estimation on throughput, but does so assuming a specific probability distribution of the channel-state and does not study the stability of the data queues either. \cite{shak2} attempts at modeling channel- and queue-state uncertainty by considering the case where only heterogeneously delayed information is available at the scheduler. They however assume knowledge of the channel-state transition probabilities. In \cite{main}, the authors study scheduling with rate adaptation in a single-hop network with a single channel under channel uncertainty. They consider cases when the channel estimates are inaccurate but complete or incomplete knowledge of the channel-estimator joint statistics is available at the scheduler. The authors, however, assume that the channel-estimates are independent across the channels for each user.

Delay performance of various wireless systems has also been investigated by many researchers recently. Among previous work in the area, the authors in \cite{prasad} study the problem of opportunistic scheduling of a wireless channel while also trying to minimize the mean delay. Neely \cite{nee3} has given a $O(1/(1-\rho))$ delay bound in the case of ON/OFF channels and a $O(N/(1-\rho))$ bound for multi-rate channels for the classical MaxWeight algorithm for a network of size $N$ and any traffic input-rate vector within a $\rho$-scaled version of the stability region (where $0 < \rho < 1$). Subsequent work in \cite{modi1} established $O(1/(1-\rho))$ and $O(N/(1-\rho))$ delay bounds for the case of single-hop and multi-hop networks, respectively, for ON/OFF channels and under both \emph{i.i.d.} and Markov modulated arrival traffic scenarios. \cite{ness3} derives lower and upper bounds on the delay in a wireless system with single-hop traffic and general interference constraints.

Our contributions of this paper are as follows:
\begin{itemize}

\item Firstly, we model the channel-estimate inaccuracy and characterize the network stability region. Compared to \cite{main}, we allow the channel estimates to have dependence among themselves, which is a more realistic situation in a modern LTE or WiMax network. Besides, we study a multi-channel setup whereas they consider a single channel.

\item Secondly, we propose two simple MaxWeight based scheduling schemes that achieve any rate in the interior of the stability region. 

\item Thirdly, we derive an $O(N/(1-\rho))$ delay bound for our system under one of the throughput-optimal policies we propose.

\item Lastly, we propose a throughput-optimal policy for the network under traffic with heterogeneous Quality of Service (QoS) constraints and present some numerical results studying its performance.

\end{itemize}


The remainder of the paper is organized as follows. In Section II, we describe the system model and the assumptions made on the arrival and channel processes. Section III provides the network stability region. We also discuss an example that illustrates that partial channel information may lead to a loss of throughput. In Section IV, we propose two throughput-optimal policies for the network and prove their optimality. We also present some simulation results in this section studying their performance. Section V provides a bound on the mean packet delay in the network. Section VI gives a throughput-optimal policy for the network under QoS constraints and studies its performance. We conclude with Section VII.

\section{System Model}

Consider a multi-user cellular downlink system with $N$ users and $M$ \emph{orthogonal} channels. The system operates with fixed-size data packets and in synchronized time slots denoted by $t \in \{0, 1, 2, \ldots\}$. It may, for example, be a cellular downlink Orthogonal Frequency-Division Multiple Access (OFDMA) system. Each user has a separate queue for its data and $Q_{i}(t)$ denotes the queue-length (in terms of packets) for user $i$ in slot $t$ where $i \in \{1,\ldots,N\}$. We assume an infinite buffer at each queue. $A_{i}(t)$ denotes the number of exogenous packet arrivals for user $i$ in slot $t$. 
It is assumed that $\{A_{i}(t)\}$ is \emph{i.i.d.} from slot to slot with $\mathbb{E}\{A_{i}(0)\}=\lambda_{i}$ and with $\mathbb{E}\{A_{i}(0)^{2}\} < \infty$ for all $i$. However, in a particular slot $t$, $A_{i}(t)$ may be dependent among themselves. 

The channel-state for user $i$ and channel $j$ in slot $t$ is denoted by $X_{ij}(t)$. We assume that $X_{ij}(t)$ is \emph{i.i.d.} from slot to slot and independent across users. Such an assumption holds, for example, in a wireless system like LTE where our \emph{channel} corresponds to a PRB in the LTE system. A PRB has 180 kHz bandwidth which is close to the coherence bandwidth of the channel for a typical delay spread of 4-5 $\mu$s (Sec. 5.3.2 in \cite{lte}). We define channel-state as the maximum number of packets that can be sent over the channel successfully without suffering an outage. We assume that $X_{ij}(t) \in \mathcal{X}$ where $\mathcal{X}$ is a discrete state-space and $\mathbb{E}\{X_{ij}^{2}(t)\} < \infty$. The scheduler has access to only estimates $S_{ij}(t)$ of the channel-state in slot $t$ where $i \in \{1,\ldots,N\}$ and $j \in \{1,\ldots,M\}$. These estimates are used by the scheduler to schedule different channels to the users in slot $t$. The estimates for a particular user may be dependent on each other. This can be used to model the effect of \emph{averaging} (or even calculating any deterministic function, for that matter) of the channel-gains as done in LTE systems which are sent by the users to the BS to reduce the feedback traffic. The only constraint we impose on $S_{ij}(t)$ is that $S_{ij}(t) \in \mathcal{S}$, where $\mathcal{S}$ is a discrete set, and that it is \emph{i.i.d.} from slot to slot. As a shorthand, we shall use $\boldsymbol{Q}(t)$, $\boldsymbol{A}(t)$ and $\mathbf{S}(t)$ to denote the queue-length vector, arrival vector and channel-estimate matrix, respectively, at slot $t$. We use $\mathbb{P}_{\mathbf{S}}(\cdot)$ to denote the probability mass function of the random variable $\mathbf{S}$. We also assume that the channel/estimator statistics given by the set of probabilities $\mathbb{P}(X_{ij}=x\ |\ \mathbf{S}=\mathbf{s})$, $\forall x \in \mathcal{X}$ and $\mathbf{s} \in \mathcal{S}$ is available at the scheduler. This can be achieved, possibly, using a mechanism that learns the statistics \emph{on-the-fly}. We assume that a channel can be allocated to at most one user in a particular slot. We use the notation $I_{j}(t)$ and $R_{j}(t)$, $j \in \{1,\ldots,M\}$ to denote the user scheduled on channel $j$ and the corresponding rate allocated to it, respectively, in the slot $t$.

We can then write the queue evolution equation as:
\begin{equation}
\label{Q_evolution}
Q_{i}(t+1) = (Q_{i}(t) - \mu_{i}(t))^{+} + A_{i}(t)
\end{equation}
where $a^{+}=max\{a,0\}$ and
\begin{equation}
\label{mu_def}
\mu_{i}(t) \triangleq \sum_{j=1}^{M} \mathbf{1}(I_{j}(t)=i) R_{j}(t) \mathbf{1}(R_{j}(t) \leqslant X_{ij}(t)).
\end{equation}
In this equation, we have assumed that the packets sent on channel $j$ are received successfully if and only if $R_{j}(t) \leqslant X_{ij}(t)$, \emph{i.e.}, probability of error is \emph{negligible}\footnote{To be precise, we can transmit data at rate $X_{ij}(t)$ with any arbitrarily small probability of error (provided $X_{ij}(t)$ is less than the Shannon capacity) assuming the physical-layer coding scheme supports it. For example, we can use an appropriate Turbo code or LDPC code for our purpose.} if $i$ transmits at a rate less than or equal to $X_{ij}(t)$ on channel $j$. Under these conditions, $\{\boldsymbol{Q}(t)\}$ is a countable state Markov chain. For simplicity, we will assume it to be irreducible. But the general case can be easily handled. 

We will use the following notation: $\mathcal{C}[A]$ is the convex hull and $int[A]$ is the interior of set $A$, $\vec{\boldsymbol{1}}_{i}$ is the $i^{\mathrm{th}}$ coordinate vector and $\vec{\boldsymbol{0}}$ and $\vec{\boldsymbol{1}}$ denote $N$-dimensional vectors of zeroes and ones, respectively.

Let $L(\boldsymbol{Q}(t))$ be a Lyapunov function. We define one-slot conditional Lyapunov drift as
\begin{equation}
\label{delta}
\Delta(\boldsymbol{Q}(t)) \triangleq \mathbb{E}\{L(\boldsymbol{Q}(t + 1)) - L(\boldsymbol{Q}(t))\ |\ \boldsymbol{Q}(t)\}.
\end{equation}
In the following, we provide an upper bound on $\Delta(\boldsymbol{Q}(t))$ which will be used later on.

\begin{lem}
For the quadratic Lyapunov function $L(\boldsymbol{Q}(t)) = \sum_{i=1}^{N} Q_{i}^{2}(t)$,
\begin{align}
\label{delta0}
\Delta(\boldsymbol{Q}(t)) & \leqslant \mathbb{E}\{U(t)\ |\ \boldsymbol{Q}(t)\} + 2 \sum_{i=1}^{N} Q_{i}(t)(\lambda_{i} - \mathbb{E}\{\mu_{i}(t)\ |\ \boldsymbol{Q}(t)\})
\end{align}
where
\begin{equation}
\label{U_t_def}
U(t) \triangleq \sum_{i=1}^{N} (A_{i}^{2}(t) + \mu_{i}^{2}(t)).
\end{equation}
Further, if $\mathbb{E}\{A_{i}^{2}(t)\} < \infty$ and $\mathbb{E}\{X_{ij}^{2}(t)\} < \infty$ for each $i$ and $j$, there exists $B < \infty$ such that
\begin{equation}
\label{delta1}
\Delta(\boldsymbol{Q}(t)) \leqslant B + 2 \sum_{i=1}^{N} Q_{i}(t)(\lambda_{i} - \mathbb{E}\{\mu_{i}(t)\ |\ \boldsymbol{Q}(t)\}).
\end{equation}
\end{lem}

\section{The Network Stability Region}
We first define the notion of stability we use in the paper.
\begin{defn}
A queue $Q_{i}(t)$ is \emph{strongly stable} if
\begin{equation*}
\limsup_{t \rightarrow \infty} \frac{1}{t} \sum_{\tau=0}^{t-1} \mathbb{E} \left\lbrace Q_{i}(\tau) \right\rbrace < \infty.
\end{equation*}
The network of queues is \emph{strongly stable} if each individual queue is strongly stable.
\end{defn}
\emph{Strong stability} implies positive recurrence of the Markov chain $\{\boldsymbol{Q}(t)\}$. In general, it is a stronger notion than positive recurrence. Throughout the paper, we shall use the term ``stability" to refer to \emph{strong stability}.

We characterize the network stability region of the system now. Consider the set of stationary policies G that base their scheduling decisions at time $t$ only on $(\boldsymbol{Q}(t)$, $\mathbf{S}(t))$ and the channel-estimator statistics. The network stability region is defined to be the closure of the arrival rates that can be stably supported by the policies in G.
\begin{thm}
The network stability region $\Lambda$ is given by
\begin{align*}
\Lambda & = \sum_{\mathbf{s} \in \mathcal{S}} \mathbb{P}_{\mathbf{S}}(\mathbf{s})\ \mathcal{C} \Bigg{[}\vec{\ \boldsymbol{0}}, \sum_{j=1}^{M}\ \mathbb{P}(X_{ij} \geqslant r_{ij}^{*}(\mathbf{s})\ |\ \mathbf{S}=\mathbf{s})\ r_{ij}^{*}(\mathbf{s})\ \vec{\boldsymbol{1}}_{i}; \forall\ i \in \{1,\ldots,N\} \Bigg{]},
\end{align*}
where
\[
r_{ij}^{*}(\mathbf{s}) \triangleq \mathop{\arg\max}\limits_{x \in \mathcal{X}}\ \{\ \mathbb{P}(X_{ij} \geqslant x\ |\ \mathbf{S}=\mathbf{s})\ x \},
\]
ties being broken lexicographically.
\end{thm}

\textbf{Proof}
The proof goes along the lines of the proof of Proposition 1 in \cite{main}. We have included it in Appendix A for the sake of completeness.

\subsection{An example}
We illustrate the loss in throughput caused by partial channel information using a simple example. We show that scheduling schemes that naively trust the channel-estimate fed back as the \emph{true} channel-state may perform much worse than the policies that don't.

Consider a system with $1$ user and $2$ channels. The $2$ channel are assumed independent. Also, $\mathbb{P} \{ X_{11}(t) = 0 \} = 0.5 = \mathbb{P} \{ X_{11}(t) = 2 \}$ and $\mathbb{P} \{ X_{12}(t) = 0 \} = \mathbb{P} \{ X_{12}(t) = 6 \} = 0.5$. Suppose we can only observe the arithmetic average of the two channel states and not the individual states. So, $\mathbb{P} \{ S_{11}(t) = s \} = 0.25 = \mathbb{P} \{ S_{12}(t) = s \}$, for $s \in \{0, 1, 3, 4\}$. Now, if we take the $\mathbf{S}$ value to be the true channel-state, it can be easily shown that the mean service provided will be $\frac{1}{4} \times (0 + 1 + 3 + 4) = 2$ packets per slot. However, due to the special choice of the support set, the $\mathbf{S}$ values give us complete information about the channel-state of both the channels. Then, it is easy to see that we can provide mean service of $\frac{1}{4} \times (0 + 2 + 6 + 8) = 4$ packets per slot.

We note that a careful scheduling decision can even double the mean service rate as shown in the example. Though the example may appear a bit contrived, numerical studies in the next section show that performance gains due to \emph{clever} scheduling may indeed be substantial in many realistic situations.

\section{Throughput-optimal policies}
In this section, we describe two throughput-optimal policies and also prove their optimality. Even though the \emph{STAT} policy described in the proof of Theorem 1 is throughput-optimal, we require knowledge of the arrival rates $\boldsymbol{\lambda}$ for \textit{STAT} to be able to perform the channel-allocation. The throughput-optimal policies described here, in contrast, just require the arrival rate vector to \emph{lie within} the stability region (without knowing $\boldsymbol{\lambda}$) and need only knowledge of the current queue-lengths. This will be available to a downlink scheduler used at a BS. Both, of course, require knowledge of the channel-estimator statistics. Moreover, as will be seen later in the section, scheduling schemes that naively trust the channel-estimate fed back perform worse than the policies in this section. For notational simplicity, we shall drop all the slot indices in this section. Lyapunov drift analysis techniques are used to prove the throughput-optimality of the policies in this section.

\subsection{MaxWeight policy}
We consider the MaxWeight (MW) type policy described below.

At each slot $t$, the channel-estimate $\mathbf{S} = \mathbf{s}$ is observed, and the decisions $I_{j}$ and $R_{j}$ are computed separately, for each $j$, as follows:
\begin{enumerate}

\item To each user $i$, assign rate $R_{ij}$ such that,
\[
R_{ij} = \mathop{\arg\max}\limits_{x \in \mathcal{X}}\ \{\ \mathbb{P}(X_{ij} \geqslant x\ |\ \mathbf{S}=\mathbf{s})\ x \}.
\]

\item Schedule the user $I_{j}$ that maximizes the \emph{rate-backlog-success-probability product}:
\[
I_{j} = \mathop{\arg\max}\limits_{1 \leqslant i \leqslant N}\ \{Q_{i}\ \mathbb{P}(X_{ij} \geqslant R_{ij}\ |\ \mathbf{S}=\mathbf{s})\ R_{ij} \}.
\] 

\end{enumerate}
For the sake of completeness, we assume that all ties here are broken lexicographically.

\begin{thm}
The MaxWeight policy is throughput-optimal.
\end{thm}

\textbf{Proof}
See Appendix B.

\subsection{Iterative MaxWeight policy}
We now analyses an iterative version of the above MaxWeight policy. These policies have also been studied in \cite{SSG} and \cite{shak1}. The new policy will be referred to as \emph{iMW}. We study this policy because we find that it can give a lower mean delay than \emph{MW} in some networks. 
In the \textit{iMW} policy, we allocate the channels sequentially from $1$ to $M$ in $M$ rounds taking into account the channels allocated so far. The \emph{virtual} queue-lengths at the beginning of round $j$ are considered for the allocation in round $j$. To aid the analysis, we use $Q_{i}^{(j)}$ to denote the \emph{virtual} queue-length of queue $i$ at the beginning of round $j$ of the allocation. $Q_{i}^{(1)}$ is defined to be $Q_{i}$. We also assume that the set $\mathcal{X}$ contains a largest element denoted by $x_{max}$. We can formulate the \emph{iMW} policy now as follows.

At each slot $t$, the channel-estimate $\mathbf{S} = \mathbf{s}$ is observed. The decisions $I_{j}$ and $R_{j}$ are computed sequentially from channel $1$ to $M$, as follows:
\begin{enumerate}

\item Start with $j = 1$.

\item For each $j$, do the following:
\begin{enumerate}

\item To each user $i$, assign rate $R_{ij}$ such that,
\[
R_{ij} = \mathop{\arg\max}\limits_{x \in \mathcal{X}}\ \{\ \mathbb{P}(X_{ij} \geqslant x\ |\ \mathbf{S}=\mathbf{s})\ x \}.
\]

\item Schedule the user $I_{j}$ that maximizes the \emph{rate-backlog-success-probability product}:
\[
I_{j} = \mathop{\arg\max}\limits_{1 \leqslant i \leqslant N}\ \{Q_{i}^{(j)}\ \mathbb{P}(X_{ij} \geqslant R_{ij}\ |\ \mathbf{S}=\mathbf{s})\ R_{ij} \}.
\] 

\end{enumerate}

\item If $j = M$, stop.\\
Else, put $Q_{i}^{(j+1)} = (Q_{i}^{(j)} - \mathbf{1}(I_{j}=i) R_{j})^{+}$, for $1 \leqslant i \leqslant N$, increment $j$ by $1$ and continue.

\end{enumerate}
Here, as before, ties are broken lexicographically.


\begin{thm}
The iterative MaxWeight policy is throughput-optimal.
\end{thm}

\textbf{Proof}
See Appendix C.

\subsection{Simulations}

\begin{figure}[!htbp]
\centering
\includegraphics[trim=1mm 3mm 9mm 6mm, clip=true, scale=0.7]{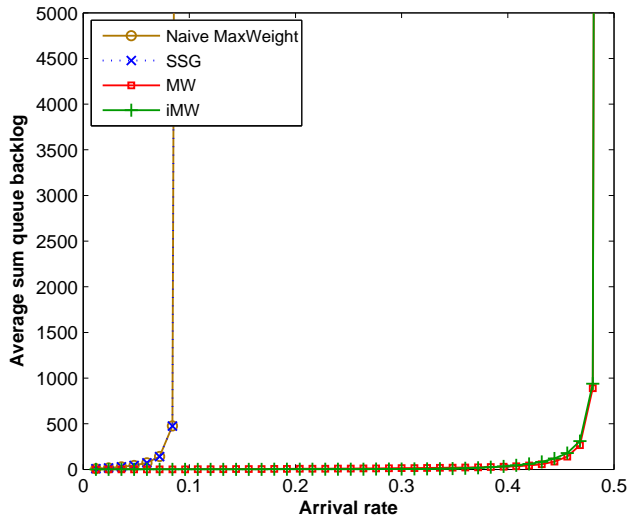}
\caption{Simulation for the ON/OFF system with symmetric traffic}
\label{fig_first_case}
\end{figure}

We present numerical simulation results to compare \emph{MW} with \emph{iMW} and also show the advantage of using channel estimators instead of the average channel gains. Firstly, we consider an \emph{ON/OFF} system with $\mathbb{P} \{ X_{ij}(t) = 0 \} = 1/2 = \mathbb{P} \{ X_{ij}(t) = 1 \}$, for all $i \in \{1,\ldots,10\}$ and $j \in \{1,\ldots,6\}$. We assume that each user estimates its channel-state correctly all the time but feeds back only the sum of the six channel-estimates to the scheduler. We do this to study the effect of \emph{averaging} the channel-estimates as in LTE\footnote{We note that this is somewhat different from an LTE setup. Arithmetic mean or EESM \cite{lte} is usually used in that case.} on the stability region. The \emph{naive} scheduling schemes (MaxWeight\cite{sri3} and SSG\cite{SSG}) calculate the average channel-estimate from the sum fed back, round it down to nearest integer and use that for scheduling, taking it to be the \emph{true} channel-state for each of the six channels. We consider symmetric Binomial$(10, \lambda)$ arrivals with equal rates for all users. We have simulated the system for $10^6$ slots for values of $\lambda$ from $0.01$ to $0.5$. The resulting simulated queue backlogs are shown in Fig. \ref{fig_first_case}. We see a huge gain in the stability region compared to the \emph{naive} algorithms. Similar results are obtained when the channel-estimates are rounded up to nearest integer instead of down.

\begin{figure}[!htbp]
\centering
	\includegraphics[trim=2mm 3mm 9mm 6mm, clip=true, scale=0.7]{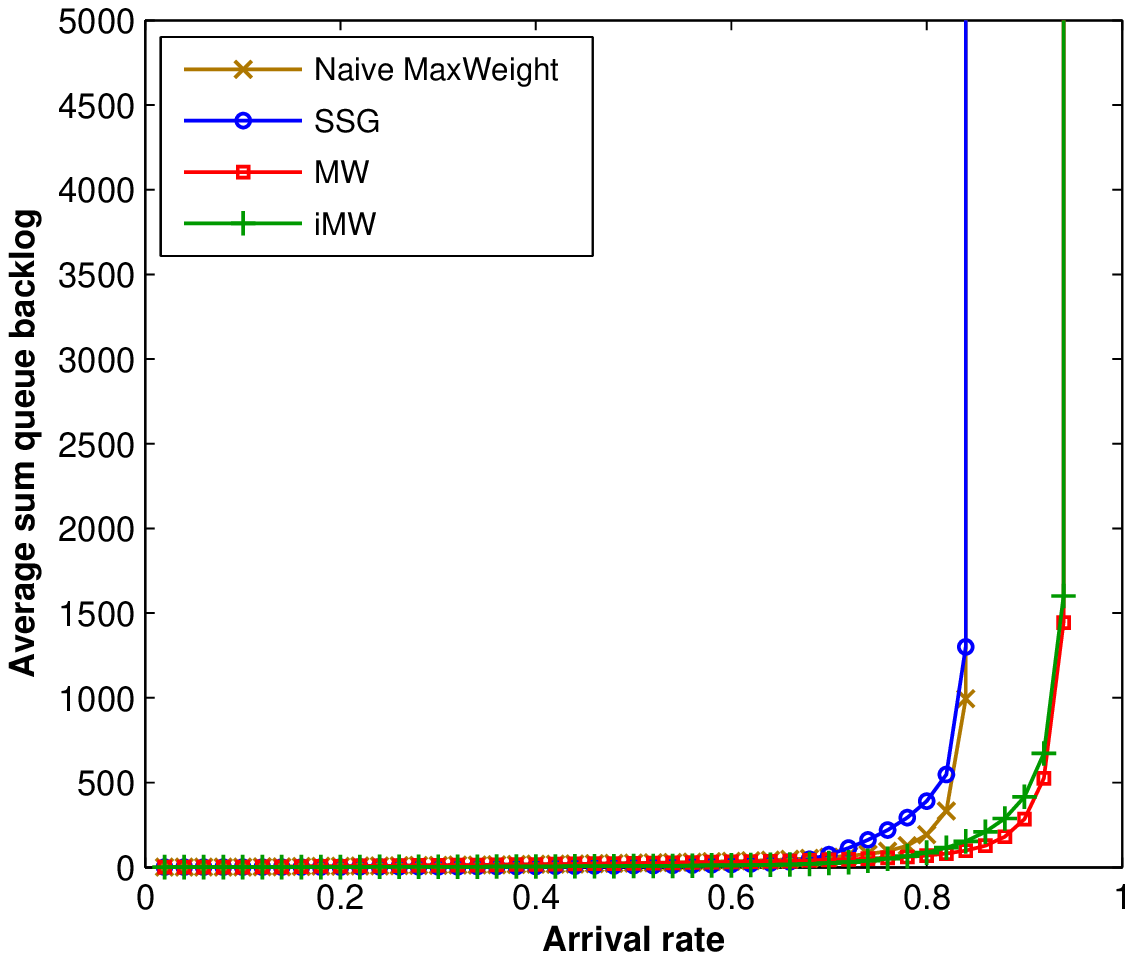}
    \llap{
    		\raisebox{12mm}{
      		\put(-200,0){
      			\includegraphics[trim=6.5mm 2mm 8mm 4mm, clip=true, scale=0.7]{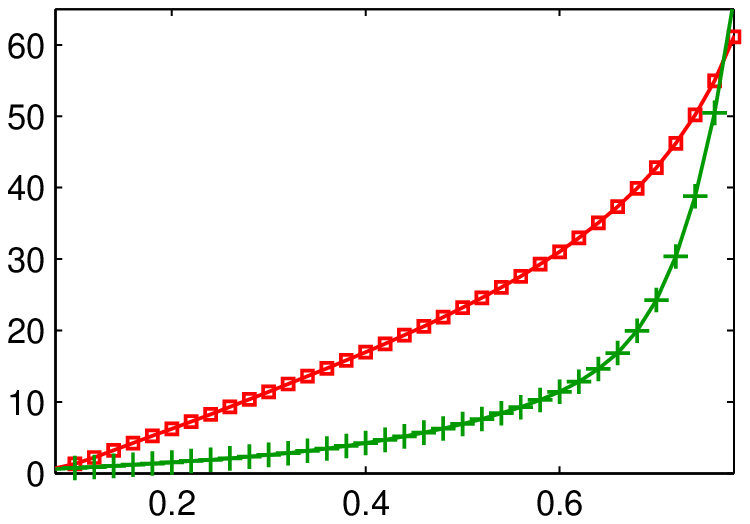}
      		}
    		}
    	}
\caption{Simulation for the multi-rate system with symmetric traffic. The inset shows the queue backlog performance of \emph{iMW} versus \emph{MW} under moderate traffic conditions.}
\label{fig_second_case}
\end{figure}

In the second set of simulations, we consider a multi-rate system with $\mathbb{P} \{ X_{ij}(t) = x \} = 1/4$, for all $i \in \{1,\ldots,10\}$, $j \in \{1,\ldots,6\}$ and $x \in \{0,1,2,3\}$. Here too, we study the \emph{averaging} effect by making the \emph{naive} scheduling schemes (MaxWeight and SSG) calculate the average channel-estimate from the sum fed back, round it up to the nearest integer and use that for scheduling. We simulate the system for $10^7$ slots with symmetric Binomial$(10, \lambda)$ arrivals for values of $\lambda$ from $0.05$ to $1$. Fig. \ref{fig_second_case} shows the resulting mean queue backlogs. We again see a decrease in the stability region when we use the \emph{naive} schedulers. In this figure, we have also expanded the graph to show that \emph{iMW} performs better than \emph{MW} at least under moderate traffic conditions. Similar results are obtained when the channel-estimates are rounded down instead of up.

We also simulated asymmetric systems and a similar behavior of the corresponding stability regions was observed. The results are not reported here for lack of space.

\section{A delay bound}
In this section, we derive an upper bound on the mean delay for the \emph{MW} policy. 

\begin{thm}
Assume $\boldsymbol{\lambda} \in int[\Lambda]$. Then, the average delay in the system, denoted by $\overline{D}$, satisfies the following bound:
\begin{equation}
\label{delay_bound}
\overline{D} \leqslant \frac{N (1 + K) \sum_{i=1}^{N} \mathbb{E}\{A_{i}^{2}(\tau)\}}{2 \mu (1 - \rho) \sum_{j=1}^{N} \mathbb{E}\{A_{j}(\tau)\}}
\end{equation}
where 
\begin{equation}
\label{def_mu}
\mu \triangleq \min_{1 \leqslant i \leqslant N} \left( \sum_{\mathbf{s} \in \mathcal{S}} \mathbb{P}_{\mathbf{S}}(\mathbf{s}) \sum_{j=1}^{M} \mathbb{P}(X_{ij} \geqslant r_{ij}^{*}(\mathbf{s})\ |\ \mathbf{S}=\mathbf{s}) r_{ij}^{*}(\mathbf{s}) \right),
\end{equation}
and $0 < \rho < 1$ and $K < \infty$.
\end{thm}

\textbf{Proof}
See Appendix D.

Note that the delay bound derived is $O(N/(1 - \rho))$ since $\frac{\sum_{i=1}^{N} \mathbb{E}\{A_{i}^{2}(\tau)\}}{\sum_{j=1}^{N} \mathbb{E}\{A_{j}(\tau)\}} = O(1)$, and $K$ and $\mu$ are fixed constants.

\section{QoS constraints}
In the previous sections, there were no QoS constraints on the data. We only had to ensure that all the queues were stable. In this section, we consider three types of traffic:
\begin{itemize}
\item Real-time (RT) traffic: Demands an upper bound on packet dropping ratio and delay deadline guarantees.
\item Rate-guarantee (RG) traffic: Demands minimum rate guarantees.
\item Best-effort (BE) traffic: Demands only queue stability.
\end{itemize}

Let $\mathcal{R}$, $\mathcal{G}$  and $\mathcal{B}$ denote the set of RT, RG and BE users, respectively, in the system. The set of all users is denoted by $\mathcal{N}$. We have a slotted system as before. However, the slots are grouped into frames and each frame consists of $T$ consecutive slots. We assume that the RT packets have a deadline of $T$ slots. These RT packets might be coming, for example, from a voice or a video source. RG packets coming into the system only require guarantees on minimum rate. These packets might be coming from a source with flow-control, such as by TCP protocol. Thus, the RG users may be treated, for our purposes, as packet sources with infinite backlogs. The minimum-rate guarantee sought by a RG user $g$ is denoted by $\beta_g$. The channel model is the same as before except that we assume $X_{bj}(t) \leqslant x_{max} < \infty$, $\forall b \in \mathcal{B}$ and for each channel $j$. All exogenous packets arriving into the system arrive \emph{only} at the start of the frame and $A_n[k]$ denotes the number of packet arrivals in frame $k$, $k \in \{0,1,2,\ldots\}$ for user $n \in \mathcal{N}$. As before, $\{A_{n}[k]\}$ is \emph{i.i.d.} from frame to frame with $\mathbb{E}\{A_{n}[0]\}=\lambda_{n}$ and with $\mathbb{E}\{A_{n}[0]^{2}\} < \infty$ for all $n \in \mathcal{N}$. In case, some of the RT packets arriving in a frame could not be served by the end of the frame, they are simply dropped. We denote the packet dropping ratio for user $r$ by $\alpha_r$, \emph{i.e.}, at most $\alpha_r$ fraction of the packets arriving for user $r$ can be dropped in the long term.

In order to satisfy the QoS constraints of the RT and RG users, we use the concept of \emph{virtual queues} (see \cite{lemma}, \cite{neely_drpc}) that evolve from frame to frame. Corresponding to each RT user $r$ and RG user $g$, we have virtual queue-length processes $\{Y_r[k]\}_{k=0}^{\infty}$  and $\{Z_g[k]\}_{k=0}^{\infty}$, respectively. We can then write the queue evolution equations for the RT users as,
\begin{equation}
\label{Y_update}
Y_r[k+1] = (Y_r[k] - \mu_r[k] + A_r[k](1-\alpha_r))^{+}, \forall r \in \mathcal{R},
\end{equation}
where $Y_r[0]=0$. $\mu_r[k] \triangleq \sum_{t=kT}^{kT+T-1} \mu_r(t)$, and $\mu_r(t)$ is as defined in \eqref{mu_def}. Similarly, the virtual queues of the RG users are updated as follows:
\begin{equation}
\label{Z_update}
Z_g[k+1] = (Z_g[k] - \mu_g[k] + \beta_g)^{+}
\end{equation}
$\forall g \in \mathcal{G}$, where $Z_g[0]=0$ and $\mu_g[k]$ is defined just like $\mu_r[k]$. Unlike the RT and RG users, the queues of the BE users evolve from slot to slot. Besides, BE users only maintain real queues since they do not have any QoS constraints whatsoever. $Q_b(t)$ denotes the queue-length of BE user $b$ at slot $t$ and its dynamics are governed by the equations,
\begin{subequations}
\label{BE_Q_evolution}
\begin{align}
Q_b(kT+1) &= (Q_b(kT) - \mu_b(kT) + A_b[k])^{+},\\
Q_b(kT+t+1) &= (Q_b(kT+t) - \mu_b(kT+t))^{+},
\end{align}
\end{subequations}
$\forall t \in \{1, 2, \ldots, T-1\}$ in frame $k$. We also define:
\begin{equation}
\label{def}
\Phi_n[k] =
\begin{cases} 
Y_n[k], & \text{if $n \in \mathcal{R}$,}\\
Z_n[k], & \text{if $n \in \mathcal{G}$, and}\\
Q_n(kT), & \text{if $n \in \mathcal{B}$.}
\end{cases}
\end{equation}


\subsection{Throughput-optimal policy}
We consider a modified version of the MaxWeight type policy described in the previous sections, which we call QoS-MaxWeight (QMW).

At each slot $t$ in frame $k$, the channel-estimate $\mathbf{S}(t) = \mathbf{s}$ is observed, and the decisions $I_{j}(t)$ and $R_{j}(t)$ are computed separately for each channel as follows:
\begin{enumerate}

\item To each user $n$ and channel $j$, assign rate $R_{nj}$ such that,
\[
R_{nj} = \mathop{\arg\max}\limits_{x \in \mathcal{X}}\ \{\ \mathbb{P}(X_{nj}(t) \geqslant x | \mathbf{S}(t)=\mathbf{s})\ x \}.
\]

\item Schedule the user $I_{j}(t)$, on channel $j$, that maximizes the \emph{rate-backlog-success-probability product} below:
\[
I_{j}(t) = \mathop{\arg\max} \limits_{n \in \mathcal{N}} \{ \Phi_n[k] \mathbb{P}(X_{nj}(t) \geqslant R_{nj} | \mathbf{S}(t)=\mathbf{s}) R_{nj} \}.
\] 

\end{enumerate}
For the sake of completeness, we assume that all ties here are broken lexicographically.

\begin{thm}
The QoS-MaxWeight policy is throughput-optimal.
\end{thm}

\textbf{Proof}
See Appendix E.

\subsection{Simulation results}

\begin{figure}[!htbp]
\centering
\includegraphics[trim=10mm 3mm 11mm 3mm, clip=true, scale=0.65]{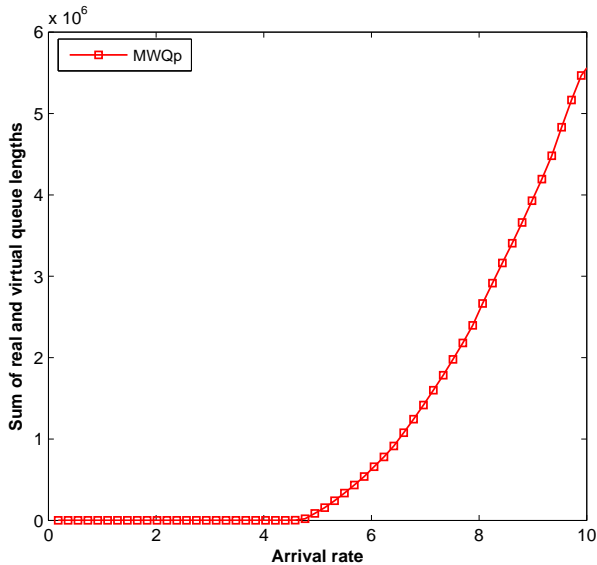}
\caption{Simulation for the multi-rate system with partial channel information and with RT, RG and BE traffic.}
\label{fig_mqos}
\end{figure}

We investigate the effect of partial information on the network stability region using numerical simulation. We consider a multi-rate system with $\mathbb{P} \{ X_{ij}(t) = x \} = 1/4$, for all $i \in \{1,\ldots,15\}$, $j \in \{1,\ldots,6\}$ and $x \in \{0,1,2,3\}$. There are $2$ RT users (with the maximum dropping ratios and deadlines being $0.01$ and $10$ slots for user $1$, and $0.02$ and $10$ slots for user $2$), $3$ RG users (with minimum rate guarantees being $5$, $2$ and $1$ packets-per-frame for users $3$, $4$ and $5$, respectively) and the rest $10$ are BE users. RT users have Binomial$(10$, $2.75)$ arrivals. RG users have Binomial arrivals with parameters $10$ and the corresponding minimum rate. We have simulated the system for $5 \times 10^6$ slots with symmetric Binomial$(10, \lambda)$ arrivals for the BE users with the value of $\lambda$ varying from $0.18$ to $10$. Fig. \ref{fig_mqos} shows the resulting mean sum queue backlogs.

\section{Conclusions}
We have considered a multi-rate wireless downlink with multiple channels and users. The base station schedules the traffic based on partial channel information. We have obtained its network stability region. We have then proposed two throughput-optimal schemes to achieve the stability region and proved their optimality. We have also derived a bound on the mean packet delay in the network. Finally, we have proposed a throughput-optimal policy for the network under QoS constraints. A natural extension of our work would be to investigate the effect of uncertainty in the queue-length information. Also, more complicated traffic models are left to be studied in greater detail.




\section*{Appendix A}

\textbf{Proof of Theorem 1:}

\textit{Sufficiency}: We show that $\boldsymbol{\lambda} \in int[\Lambda]$ is a sufficient condition for stability. Now, since $\Lambda$ is convex, for each $\mathbf{s} \in \mathcal{S}$, there exists a scaling vector $\boldsymbol{\gamma}^{\mathbf{s}}$ and a scalar $\epsilon > 0$ such that
\begin{equation}
\label{conv}
\lambda_{i} + \epsilon < \sum_{\mathbf{s} \in \mathcal{S}} \gamma^{\mathbf{s}}_{i}\ \mathbb{P}_{\mathbf{S}}(\mathbf{s}) \sum_{j=1}^{M} \mathbb{P}(X_{ij} \geqslant r_{ij}^{*}(\mathbf{s})\ |\ \mathbf{S}=\mathbf{s})\ r_{ij}^{*}(\mathbf{s})
\end{equation}
for any user $i$, and where $\sum_{i=1}^{N} \gamma_{i}^{\mathbf{s}} = 1,\forall \mathbf{s} \in \mathcal{S}$.

Consider the following stationary randomized policy, henceforth referred to as \emph{STAT}: for channel estimate $\mathbf{S}(t) = \mathbf{s}$, allocate all channels to user $i$ with probability $\gamma_{i}^{\mathbf{s}}$ and set the rate allocated on channel $j$ to it to be $r_{ij}^{*}(\mathbf{s})$. In that case, the service rate of user $i$, call it $u_{i}$, is given by
\begin{align}
\label{stat}
u_{i} & = \mathbb{E}\{\mu_{i}(t)\} \nonumber \\
& = \mathbb{E}\{\mu_{i}(t)\ |\ \boldsymbol{Q}(t)\} \nonumber \\
& = \sum_{\mathbf{s} \in \mathcal{S}} \gamma^{\mathbf{s}}_{i}\ \mathbb{P}_{\mathbf{S}}(\mathbf{s}) \sum_{j=1}^{M} \mathbb{P}(X_{ij} \geqslant r_{ij}^{*}(\mathbf{s})\ |\ \mathbf{S}=\mathbf{s})\ r_{ij}^{*}(\mathbf{s}).
\end{align}

Considering the quadratic Lyapunov function defined before and using \eqref{conv} and \eqref{stat} in \eqref{delta1}, we can write the Lyapunov drift inequality as
\begin{align}
\label{deldel}
\Delta(\boldsymbol{Q}(t)) & < B - 2 \epsilon \sum_{i=1}^{N} Q_{i}(t).
\end{align}
Now, we note that $\boldsymbol{Q}(t)$ evolves as a Markov chain with $L(\boldsymbol{Q}(t))$ bounded below by zero and $\mathbb{E}\{L(\boldsymbol{Q}(t+1)) \ |\ \boldsymbol{Q}(t)\} < \infty$. Also, notice that the conditional drift in \eqref{deldel} is negative outside the finite set $\{ \boldsymbol{Q}(t) : \sum_{i=1}^{N} Q_{i}(t) \leqslant \frac{B}{2 \epsilon}\}$. Thus, using Theorem 2 from \cite{definition}, we can say that the queue-length process is stable.



\textit{Necessity}: Here we show that $\boldsymbol{\lambda} \in \Lambda$ is a necessary condition for stability. For if not so, i.e., if $\boldsymbol{\lambda} \notin \Lambda$, we have a vector $\boldsymbol{\alpha}$ and a scalar $\delta > 0$, such that for any $\boldsymbol{v} \in \Lambda$, we have
\[
\sum_{i=1}^{N} \alpha_{i}(\lambda_{i}-v_{i}) \geqslant \delta,
\]
from the Strict Separation Theorem (see, for example, Proposition B.14 in \cite{book}). Let us define the linear Lyapunov function $L(\boldsymbol{Q}(t)) = \sum_{i=1}^{N} \alpha_{i} Q_{i}(t)$. Now, using definition (\ref{delta}) and the queue evolution equation (\ref{Q_evolution}), for any stationary policy in G, we can write
\begin{align}
\label{delta2}
\Delta(\boldsymbol{Q}(t)) & \geqslant \sum_{i=1}^{N} \alpha_{i}\mathbb{E}\{A_{i}(t) - \mu_{i}(t)\ |\ \boldsymbol{Q}(t)\} \nonumber \\
& = \sum_{i=1}^{N} \alpha_{i}(\lambda_{i} - \mathbb{E}\{\mu_{i}(t)\ |\ \boldsymbol{Q}(t)\}).
\end{align}
Let $u_{i} \triangleq \mathbb{E}\{\mu_{i}(t)\ |\ \boldsymbol{Q}(t)\}$. We now show that $\boldsymbol{u} \in \Lambda$.
\begin{align*}
& \quad \mathbb{E}\left[\mu_{i}(t)\ |\ \boldsymbol{Q}(t)\right]\\
& = \sum_{\mathbf{s} \in \mathcal{S}} \mathbb{P}_{\mathbf{S}}(\mathbf{s}) \mathbb{E}\left[ \mu_{i}(t)\ |\ \boldsymbol{Q}(t); \mathbf{S}(t) = \mathbf{s} \right]\\
& = \sum_{\mathbf{s} \in \mathcal{S}} \mathbb{P}_{\mathbf{S}}(\mathbf{s}) \Bigg( \sum_{j=1}^{M} \mathbb{P}(R_{j}(t) \leqslant X_{ij}(t)\ |\ \mathbf{S}(t) = \mathbf{s}) \mathbf{1}(I_{j}(t)=i) R_{j}(t) \Bigg)\\
& \leqslant \sum_{\mathbf{s} \in \mathcal{S}} \mathbb{P}_{\mathbf{S}}(\mathbf{s}) \Bigg( \sum_{j=1}^{M} \mathbb{P}(r_{ij}^{*}(\mathbf{s}) \leqslant X_{ij}(t)\ |\ \mathbf{S}(t) = \mathbf{s}) \mathbf{1}(I_{j}(t)=i) r_{ij}^{*}(\mathbf{s}) \Bigg).
\end{align*}
The second equality holds since the policy decisions $I_{j}(t)$ and $R_{j}(t)$ are completely determined by $\boldsymbol{Q}(t)$ and $\mathbf{S}(t)$ within the class of stationary policies G. The inequality at the end holds because of the way $r_{ij}^{*}(\mathbf{s})$ is defined. Thus, $\boldsymbol{u} \in \Lambda$.

Thus, from the Strict Separation Theorem and (\ref{delta2}), we get
\begin{align*}
\Delta(\boldsymbol{Q}(t)) & \geqslant \sum_{i=1}^{N} \alpha_{i}(\lambda_{i} - \mathbb{E}\{\mu_{i}(t)\ |\ \boldsymbol{Q}(t)\})\\
& = \sum_{i=1}^{N} \alpha_{i}(\lambda_{i} - u_{i})\\
& \geqslant \delta.
\end{align*}
It can be shown that $\mathbb{E}\{|L(\boldsymbol{Q}(t + 1)) - L(\boldsymbol{Q}(t))| | \boldsymbol{Q}(t)\} < \infty$. Taking the finite set as $\{0\}$ and noting that $L(\boldsymbol{Q}(t)) > 0$ for at least some $\boldsymbol{Q}(t)$, we see that the DTMC will \emph{not} be positive recurrent, and hence, not \emph{strongly stable} either. Thus, we conclude that the network will be unstable.

\section*{Appendix B}

\textbf{Proof of Theorem 2:}

Assume that we are working with $\boldsymbol{\lambda} \in int[\Lambda]$. Let the \emph{MW} policy make decisions $I_{j}$ and $R_{j}$, and the \emph{STAT} policy make decisions $I_{j}^{\prime}$ and $R_{j}^{\prime}$. We then have,
\begin{align*}
& \quad Q_{i}\ \mathbf{1}(I_{j}=i)\ R_{j}\ \mathbb{P}(X_{ij} \geqslant R_{j}\ |\ \mathbf{S}=\mathbf{s}) \geqslant Q_{i}\ \mathbf{1}(I_{j}^{\prime}=i)\ R_{j}^{\prime}\ \mathbb{P}(X_{ij} \geqslant R_{j}^{\prime}\ |\ \mathbf{S}=\mathbf{s})
\end{align*}
$\forall i \in \{1,\ldots,N\}$ and $\forall j \in \{1,\ldots,M\}$, since \emph{MW} maximizes the expression in the LHS over the class G. Thus, summing over all $i$ and $j$, we get
\begin{align*}
& \quad \sum_{i=1}^{N} \sum_{j=1}^{M} Q_{i}\ \mathbf{1}(I_{j}=i)\ R_{j}\ \mathbb{P}(X_{ij} \geqslant R_{j}\ |\ \mathbf{S}=\mathbf{s}) \nonumber \\
& \geqslant \sum_{i=1}^{N} \sum_{j=1}^{M} Q_{i}\ \mathbf{1}(I_{j}^{\prime}=i)\ R_{j}^{\prime}\ \mathbb{P}(X_{ij} \geqslant R_{j}^{\prime}\ |\ \mathbf{S}=\mathbf{s}).
\end{align*}
Therefore,
\begin{align*}
& \quad \sum_{i=1}^{N} \sum_{j=1}^{M} Q_{i}\ \mathbb{E}\{\ \mathbf{1}(I_{j}=i)\ R_{j}\ \mathbf{1}(X_{ij} \geqslant R_{j})\ |\ \boldsymbol{Q}; \mathbf{S}=\mathbf{s}\} \nonumber \\
& \geqslant \sum_{i=1}^{N} \sum_{j=1}^{M} Q_{i}\ \mathbb{E}\{\mathbf{1}(I_{j}^{\prime}=i)\ R_{j}^{\prime}\ \mathbf{1}(X_{ij} \geqslant R_{j}^{\prime})\ |\ \boldsymbol{Q}; \mathbf{S}=\mathbf{s}\}
\end{align*}
where we have used the fact that the policy decisions are completely determined by $\boldsymbol{Q}$ and $\mathbf{S}$ within the class G. Taking expectation w.r.t. $\mathbf{S}$,
\begin{align}
\label{deriv1}
& \quad \sum_{i=1}^{N} \sum_{j=1}^{M} Q_{i}\ \mathbb{E}\{\mathbf{1}(I_{j}=i)\ R_{j}\ \mathbf{1}(X_{ij} \geqslant R_{j})\ |\ \boldsymbol{Q}\} \nonumber \\
& \geqslant \sum_{i=1}^{N} \sum_{j=1}^{M} Q_{i}\ \mathbb{E}\{\mathbf{1}(I_{j}^{\prime}=i)\ R_{j}^{\prime}\ \mathbf{1}(X_{ij} \geqslant R_{j}^{\prime})\ |\ \boldsymbol{Q}\}.
\end{align}
Note that \emph{STAT} stabilizes the system since $\boldsymbol{\lambda} \in int[\Lambda]$. Thus, for some $\epsilon > 0$, we can write
\begin{align*}
\lambda_{i} + \epsilon & < \sum_{j=1}^{M} \mathbb{E}\{ \mathbf{1}(I_{j}^{\prime}=i) R_{j}^{\prime} \mathbf{1}(R_{j}^{\prime} \leqslant X_{ij}) \}\\
& = \sum_{j=1}^{M} \mathbb{E}\{ \mathbf{1}(I_{j}^{\prime}=i) R_{j}^{\prime} \mathbf{1}(R_{j}^{\prime} \leqslant X_{ij})\ |\ \boldsymbol{Q}\}
\end{align*}
where the equality holds because \emph{STAT} makes its decisions independent of the queue-lengths. Multiplying $Q_{i}$ on both sides and summing over all $i \in \{1,\ldots,N\}$, we get
\begin{align}
\label{common}
\sum_{i=1}^{N} Q_{i} (\lambda_{i} + \epsilon) & < \sum_{i=1}^{N} \sum_{j=1}^{M} Q_{i} \mathbb{E}\{ \mathbf{1}(I_{j}^{\prime}=i) R_{j}^{\prime} \mathbf{1}(R_{j}^{\prime} \leqslant X_{ij})\ |\ \boldsymbol{Q}\} \nonumber \\
& \leqslant \sum_{i=1}^{N} \sum_{j=1}^{M} Q_{i} \mathbb{E}\{ \mathbf{1}(I_{j}=i) R_{j} \mathbf{1}(R_{j} \leqslant X_{ij})\ |\ \boldsymbol{Q}\}
\end{align}
where we have used (\ref{deriv1}) in the second inequality. Therefore,
\begin{align}
\label{deriv2}
& \quad \sum_{i=1}^{N} Q_{i} \left(\lambda_{i} - \sum_{j=1}^{M} \mathbb{E}\{ \mathbf{1}(I_{j}=i) R_{j} \mathbf{1}(R_{j} \leqslant X_{ij})\ |\ \boldsymbol{Q}\}\right) < - \epsilon \sum_{i=1}^{N} Q_{i}
\end{align}
Finally, using (\ref{deriv2}) in the drift inequality (\ref{delta1}), we get
\[
\Delta(\boldsymbol{Q}(t)) < B - 2 \epsilon \sum_{i=1}^{N} Q_{i}(t).
\]
Now the queue evolves as a Markov chain since the scheduling and rate allocation decisions are taken based on the \emph{current} queue-lengths and channel-estimate. The drift inequality above gives negative drift but for a finite set of queue-lengths. Hence, using Theorem 2 from \cite{definition}, the network is stable.

\section*{Appendix C}

\textbf{Proof of Theorem 3:}

Assume that we are working with $\boldsymbol{\lambda} \in int[\Lambda]$. Let the \emph{iMW} policy make decisions $I_{j}$ and $R_{j}$, the \emph{MW} policy make decisions $I_{j}^{\prime}$ and $R_{j}^{\prime}$ and the \emph{STAT} policy make decisions $I_{j}^{\prime\prime}$ and $R_{j}^{\prime\prime}$. For ease of exposition, we use $R_{ij}$ and $R_{ij}^{\prime}$ as the \emph{intermediate} allocation variables for \emph{iMW} and \emph{MW} policies, respectively. We then have,
\begin{align*}
& \quad Q_{I_{j}}\ R_{I_{j}j}\ \mathbb{P}(X_{ij} \geqslant R_{I_{j}j}\ |\ \mathbf{S}=\mathbf{s})\\
& \geqslant Q_{I_{j}}^{(j-1)}\ R_{I_{j}j}\ \mathbb{P}(X_{ij} \geqslant R_{I_{j}j}\ |\ \mathbf{S}=\mathbf{s})\\
& \geqslant Q_{I_{j}^{\prime}}^{(j-1)}\ R_{I_{j}^{\prime}j}^{\prime}\ \mathbb{P}(X_{ij} \geqslant R_{I_{j}^{\prime}j}^{\prime}\ |\ \mathbf{S}=\mathbf{s})\\
& \geqslant (Q_{I_{j}^{\prime}} - M x_{max})\ R_{I_{j}^{\prime}j}^{\prime}\ \mathbb{P}(X_{ij} \geqslant R_{I_{j}^{\prime}j}^{\prime}\ |\ \mathbf{S}=\mathbf{s})\\
& \geqslant Q_{I_{j}^{\prime}}\ R_{I_{j}^{\prime}j}^{\prime}\ \mathbb{P}(X_{ij} \geqslant R_{I_{j}^{\prime}j}^{\prime}\ |\ \mathbf{S}=\mathbf{s}) - M x_{max}^{2}.
\end{align*}
In the first inequality, we have used the fact that the \emph{virtual} queue-lengths cannot \emph{increase} with allocation rounds. The second inequality is due to the way $iMW$ does the allocation at each round. Summing both sides over all $j \in \{1,\ldots,M\}$, we get
\begin{align*}
& \sum_{j=1}^{M} Q_{I_{j}}\ R_{I_{j}j}\ \mathbb{P}(X_{ij} \geqslant R_{I_{j}j}\ |\ \mathbf{S}=\mathbf{s}) \nonumber \\
& \geqslant \sum_{j=1}^{M} Q_{I_{j}^{\prime}}\ R_{I_{j}^{\prime}j}^{\prime}\ \mathbb{P}(X_{ij} \geqslant R_{I_{j}^{\prime}j}^{\prime}\ |\ \mathbf{S}=\mathbf{s}) - M^{2} x_{max}^{2}.
\end{align*}
Therefore,
\begin{align*}
& \quad \sum_{i=1}^{N} \sum_{j=1}^{M} Q_{i}\ \mathbf{1}(I_{j}=i)\ R_{ij}\ \mathbb{P}(X_{ij} \geqslant R_{ij}\ |\ \mathbf{S}=\mathbf{s}) \nonumber \\
& \geqslant \sum_{i=1}^{N} \sum_{j=1}^{M} Q_{i}\ \mathbf{1}(I_{j}^{\prime}=i)\ R_{ij}^{\prime}\ \mathbb{P}(X_{ij} \geqslant R_{ij}^{\prime}\ |\ \mathbf{S}=\mathbf{s}) - M^{2} x_{max}^{2},
\end{align*}
where we have used the property of indicator functions. Now, using the fact that $\mathbf{1}(I_{j}=i)\ R_{ij} = \mathbf{1}(I_{j}=i)\ R_{j}$, for all $i$ for the $iMW$ policy (and a similar thing for $MW$), and proceeding in a way similar to the previous proof, we will get
\begin{align}
\label{deriv3}
& \quad \sum_{i=1}^{N} \sum_{j=1}^{M} Q_{i}\ \mathbb{E}\{\mathbf{1}(I_{j}^{\prime}=i)\ R_{j}^{\prime}\ \mathbf{1}(X_{ij} \geqslant R_{j}^{\prime})\ |\ \boldsymbol{Q}\} \nonumber \\
& \leqslant \sum_{i=1}^{N} \sum_{j=1}^{M} Q_{i}\ \mathbb{E}\{\mathbf{1}(I_{j}=i)\ R_{j}\ \mathbf{1}(X_{ij} \geqslant R_{j})\ |\ \boldsymbol{Q}\} + M^{2} x_{max}^{2}
\end{align}
Now, using a $STAT$ policy, as in the previous proof, to stabilize $\boldsymbol{\lambda} + \epsilon \vec{\boldsymbol{1}} \in int[\Lambda]$ for some $\epsilon > 0$, and using the property of \textit{MW} for the expression obtained, as done in \eqref{common}, we get
\begin{align*}
\sum_{i=1}^{N} Q_{i} (\lambda_{i} + \epsilon) & \leqslant M^{2} x_{max}^{2} + \sum_{i=1}^{N} \sum_{j=1}^{M} Q_{i} \mathbb{E}\{ \mathbf{1}(I_{j}=i) R_{j} \mathbf{1}(R_{j} \leqslant X_{ij})\ |\ \boldsymbol{Q}\}
\end{align*}
where we have used the inequality in (\ref{deriv3}). Finally, rearranging the terms in the above inequality and using it in (\ref{delta1}), we get
\[
\Delta(\boldsymbol{Q}(t)) < (B + 2 M^{2} x_{max}^{2}) - 2 \epsilon \sum_{i=1}^{N} Q_{i}(t)
\]
Now, as before, we know that the queue evolves as a Markov chain. As before, observe that the drift is negative but for a finite set of queue-lengths. Hence, the network is stable.

\section*{Appendix D}

\textbf{Proof of Theorem 4:}

Here too we drop the slot indices for the sake of brevity. Let the \emph{MW} policy make decisions $I_{j}$ and $R_{j}$, and any other policy $G^{\prime}$ in class G make decisions $I_{j}^{\prime}$ and $R_{j}^{\prime}$. Then, we can write
\begin{align}
\label{deriv4}
& \quad \sum_{i=1}^{N} \sum_{j=1}^{M} Q_{i}\ \mathbb{E}\{\mathbf{1}(I_{j}=i)\ R_{j}\ \mathbf{1}(X_{ij} \geqslant R_{j})\ |\ \boldsymbol{Q}\} \nonumber \\
& \geqslant \sum_{i=1}^{N} \sum_{j=1}^{M} Q_{i}\ \mathbb{E}\{\mathbf{1}(I_{j}^{\prime}=i)\ R_{j}^{\prime}\ \mathbf{1}(X_{ij} \geqslant R_{j}^{\prime})\ |\ \boldsymbol{Q}\}.
\end{align}
Since $\boldsymbol{\lambda} \in int[\Lambda]$, we can say that $\boldsymbol{\lambda} \in \rho \Lambda$ for some $0 < \rho <1$. Consequently, $\frac{\boldsymbol{\lambda}}{\rho} \in \Lambda$. Furthermore, defining $\mu$ as in \eqref{def_mu}, we can say that $\mu^{\prime} \vec{\boldsymbol{1}} \in \Lambda$ where $\mu^{\prime} = \frac{\mu}{N}$. Thus, owing to the convexity of $\Lambda$, we have
\[
\boldsymbol{\lambda} + (1 - \rho) \mu^{\prime} \vec{\boldsymbol{1}} \in \Lambda.
\]
Hence, we can find a $STAT$ policy such that
\begin{equation}
\label{deriv5}
\mathbb{E}\left\lbrace \sum_{j=1}^{M} \mathbf{1}(I_{j}^{\prime}=i)\ R_{j}^{\prime}\ \mathbf{1}(X_{ij} \geqslant R_{j}^{\prime})\ |\ \boldsymbol{Q} \right\rbrace = \lambda_{i} + (1 - \rho) \mu^{\prime}.
\end{equation}
Substituting (\ref{deriv5}) into (\ref{deriv4}), and then using (\ref{delta0}), we get
\begin{equation}
\Delta(\boldsymbol{Q}(t)) \leqslant \mathbb{E}\{U(t)\ |\ \boldsymbol{Q}(t)\} - 2 \mu^{\prime} (1 - \rho) \sum_{i=1}^{N} Q_{i}(t).
\end{equation}
Finally, using Lemma 4.1 of \cite{lemma} and the previous inequality, we can write
\begin{align}
& \quad \limsup_{t \rightarrow \infty} \frac{1}{t} \sum_{\tau = 0}^{t-1} \sum_{i=1}^{N} \mathbb{E}\{Q_{i}(\tau)\} \leqslant \frac{N}{2 \mu (1 - \rho)} \limsup_{t \rightarrow \infty} \frac{1}{t} \sum_{\tau = 0}^{t-1} \mathbb{E}\{U(\tau)\}.
\end{align}

Now, we notice that since the system evolves as an ergodic Markov chain with countable state-space, the time-averages are well defined, and hence, $\lim$ can be used instead of $\limsup$. Also, note that since $\mathbb{E}\{X_{ij}^{2}(t)\} < \infty$ and $N$ is finite, we can write $\mathbb{E}\{\mu_{i}^{2}(t)\} < K \mathbb{E}\{A_{i}^{2}(t)\}$ for all $i$ and $K < \infty$. So, $\mathbb{E}\{U(\tau)\} < (1 + K) \sum_{i=1}^{N}  \mathbb{E}\{A_{i}^{2}(\tau)\} < \infty $. Using this fact and Little's Theorem, we thus get the delay bound in \eqref{delay_bound}.

\section*{Appendix E}

\textbf{Proof of Theorem 5:}

Squaring both sides of \eqref{Y_update} and proceeding as in Appendix F, for each RT user $r$, we get
\begin{equation}
\label{Y_sq}
Y_r[k+1]^2 - Y_r[k]^2 \leqslant (A_r[k](1-\alpha_r))^2 + \mu_r[k]^2 + 2 Y_r[k] (A_r[k](1-\alpha_r) - \mu_r[k]).
\end{equation}
Similarly, for each RG user $g$, from \eqref{Z_update}, we get
\begin{equation}
\label{Z_sq}
Z_g[k+1]^2 - Z_g[k]^2 \leqslant \beta_g^2 + \mu_g[k]^2 + 2 Z_g[k] (\beta_g - \mu_g[k]).
\end{equation}
Squaring both sides of \eqref{BE_Q_evolution}, for each BE user $b$, we get
\begin{equation}
\label{sq_a}
Q_b(kT+1)^2 - Q_b(kT)^2 \leqslant A_b[k]^2 + \mu_b(kT)^2 + 2 Q_b(kT) (A_b[k] - \mu_b(kT)),
\end{equation}
\begin{equation}
\label{sq_b}
Q_b(kT+t+1)^2 - Q_b(kT+t)^2 \leqslant \mu_b(kT+t)^2 - 2 Q_b(kT+t) \mu_b(kT+t),
\end{equation}
for $t=\{1,\ldots,T-1\}$. Now, observe that
\[
Q_b(kT+t) \geqslant Q_b(kT) - t \mu_{max},
\]
where a finite $\mu_{max}$ exists since we assumed that $X_{bj}(t) \leqslant x_{max} < \infty$, $\forall b \in \mathcal{B}$ and for each channel $j$. So, we can write
\begin{align*}
- Q_b(kT+t) \mu_b(kT+t) & \leqslant - Q_b(kT) \mu_b(kT+t) + t \mu_{max} \mu_b(kT+t) \\
& \leqslant - Q_b(kT) \mu_b(kT+t) + t \mu_{max}^2.
\end{align*}
Using the above fact, we can rewrite \eqref{sq_b} as
\[
Q_b(kT+t+1)^2 - Q_b(kT+t)^2 \leqslant \mu_b(kT+t)^2 + 2 t \mu_{max}^2 - 2 Q_b(kT) \mu_b(kT+t).
\]
Summing both sides of the previous inequality over $t=\{1,\ldots,T-1\}$, we get
\begin{equation}
Q_b(kT+T)^2 - Q_b(kT+1)^2 \leqslant (T^2 - 1) \mu_{max}^2 - 2 Q_b(kT) \sum_{t=1}^{T-1} \mu_b(kT+t).
\end{equation}
Using \eqref{sq_a} and the above inequality, we can write
\begin{equation}
\label{Q_sq}
Q_b(kT+T)^2 - Q_b(kT)^2 \leqslant A_b[k]^2 + T^2 \mu_{max}^2 + 2 Q_b(kT) (A_b[k] - \mu_b[k]),
\end{equation}
where we have used the fact that $\mu_b(kT)^2 \leqslant \mu_{max}^2$ and $\mu_b[k] \triangleq \sum_{t=0}^{T-1} \mu_b(kT+t)$.

Consider the quadratic Lyapunov function $L(\myp{k}) = \sum_{n \in \mathcal{N}} \Phi_n[k]^2$. We now define the one-frame conditional Lyapunov drift as
\[
\Delta(\myp{k}) \triangleq \myce{L(\myp{k+1}) - L(\myp{k})}{\myp{k}}.
\]
Using \eqref{Y_sq}, \eqref{Z_sq} and \eqref{Q_sq}, we can express the drift for the QMW policy as
\begin{align}
\label{dummy_bound1}
\Delta(\myp{k}) & \leqslant \mathbb{E} \Bigg\lbrace \sum_{r \in \mathcal{R}} \left( A_r[k]^2 (1-\alpha_r)^2 + \mu_r[k]^2 \right) + \sum_{b \in \mathcal{B}} \left( A_b[k]^2 + T^2 \mu_{max}^2 \right) \nonumber \\
& \quad + \sum_{g \in \mathcal{G}} \left( \beta_g^2 + \mu_g[k]^2 \right) \Bigg| \myp{k} \Bigg\rbrace + 2\ \mathbb{E} \Bigg\lbrace \sum_{b \in \mathcal{B}} \left( Q_b(kT) (A_b[k] - \mu_b[k]) \right) \nonumber \\
& \quad + \sum_{r \in \mathcal{R}} \left( Y_r[k] (A_r[k](1-\alpha_r) - \mu_r[k]) \right) + \sum_{g \in \mathcal{G}} \left( Z_g[k] (\beta_g - \mu_g[k]) \right) \Bigg| \myp{k} \Bigg\rbrace.
\end{align}
Now, we can find a $B < \infty$ such that
\begin{align*}
B & \geqslant \mathbb{E} \Bigg\lbrace \sum_{r \in \mathcal{R}} \left( A_r[k]^2 (1-\alpha_r)^2 + \mu_r[k]^2 \right) + \sum_{b \in \mathcal{B}} \left( A_b[k]^2 + T^2 \mu_{max}^2 \right) \\
& \quad + \sum_{g \in \mathcal{G}} \left( \beta_g^2 + \mu_g[k]^2 \right) \Bigg| \myp{k} \Bigg\rbrace
\end{align*}
since $\mathbb{E}\{A_{n}^{2}[k]\} < \infty$ and $\mathbb{E}\{X_{nj}^{2}(t)\} < \infty$ for each $n \in \mathcal{N}$ and channel $j$. Also, $\{A_{n}[k]\}$ is \emph{i.i.d.} and $\mathbb{E}\{A_{n}[0]\}=\lambda_{n}$. So, we can rewrite \eqref{dummy_bound1} as
\begin{align*}
\Delta(\myp{k}) & \leqslant B + 2\ \mathbb{E} \Bigg\lbrace \sum_{r \in \mathcal{R}} \left( Y_r[k] (\lambda_r(1-\alpha_r) - \mu_r[k]) \right) + \sum_{g \in \mathcal{G}} \left( Z_g[k] (\beta_g - \mu_g[k]) \right) \nonumber \\
& \quad + \sum_{b \in \mathcal{B}} \left( Q_b(kT) (\lambda_b - \mu_b[k]) \right) \Bigg| \myp{k} \Bigg\rbrace.
\end{align*}
Rearranging the terms and using the definition \eqref{def}, we get
\begin{align}
\label{dummy_bound2}
\Delta(\myp{k}) & \leqslant B + 2\ \mathbb{E} \Bigg\lbrace \sum_{r \in \mathcal{R}} Y_r[k] \lambda_r(1-\alpha_r) + \sum_{g \in \mathcal{G}} Z_g[k] \beta_g + \sum_{b \in \mathcal{B}}  Q_b(kT) \lambda_b \Bigg| \myp{k} \Bigg\rbrace \\
& - 2\ \mathbb{E} \Bigg\lbrace \sum_{n \in \mathcal{N}} \Phi_n[k] \mu_n[k] \Bigg| \myp{k} \Bigg\rbrace.
\end{align}
From the definition of the QMW algorithm, it can be  easily shown that
\begin{align*}
& \mathbb{E} \Bigg\lbrace \sum_{n \in \mathcal{N}} \Phi_n[k] \mu_n[k] \Bigg| \myp{k}, \mathbf{S}(kT), \ldots, \mathbf{S}(kT+T-1) \Bigg\rbrace \\
& \geqslant \mathbb{E} \Bigg\lbrace \sum_{n \in \mathcal{N}} \Phi_n[k] \mu_n^\prime[k] \Bigg| \myp{k}, \mathbf{S}(kT), \ldots, \mathbf{S}(kT+T-1) \Bigg\rbrace
\end{align*}
where $\mu_n^\prime[k]$ corresponds to any other stationary randomized policy. Taking expectation \emph{w.r.t.} the channel-estimates in the frame, we get
\begin{align*}
& \mathbb{E} \Bigg\lbrace \sum_{n \in \mathcal{N}} \Phi_n[k] \mu_n[k] \Bigg| \myp{k} \Bigg\rbrace \geqslant \mathbb{E} \Bigg\lbrace \sum_{n \in \mathcal{N}} \Phi_n[k] \mu_n^\prime[k] \Bigg| \myp{k} \Bigg\rbrace.
\end{align*}
Thus, using the above fact, we can rewrite \eqref{dummy_bound2} as
\begin{align}
\Delta(\myp{k}) & \leqslant B + 2\ \mathbb{E} \Bigg\lbrace \sum_{r \in \mathcal{R}} Y_r[k] \lambda_r(1-\alpha_r) + \sum_{g \in \mathcal{G}} Z_g[k] \beta_g + \sum_{b \in \mathcal{B}}  Q_b(kT) \lambda_b \Bigg| \myp{k} \Bigg\rbrace \\
& - 2\ \mathbb{E} \Bigg\lbrace \sum_{n \in \mathcal{N}} \Phi_n[k] \mu_n^\prime[k] \Bigg| \myp{k} \Bigg\rbrace.
\end{align}
Rearranging the terms again and using definition \eqref{def}, we then get
\begin{align*}
\Delta(\myp{k}) & \leqslant B + 2 \sum_{r \in \mathcal{R}} Y_r[k] \mathbb{E} \{ (\lambda_r(1-\alpha_r) - \mu_r^\prime[k]) | \myp{k} \} \nonumber \\
& \quad + 2 \sum_{g \in \mathcal{G}} Z_g[k] \mathbb{E} \{ (\beta_g - \mu_g^\prime[k]) | \myp{k} \} + 2 \sum_{b \in \mathcal{B}} Q_b(kT) \mathbb{E} \{ (\lambda_b - \mu_b^\prime[k]) | \myp{k} \}.
\end{align*}
Now, if the system is stabilizable, \emph{i.e.} if $\boldsymbol{\lambda} \in int[\Lambda]$, then a stationary randomized policy exists that stabilizes the system. Besides, this policy makes decisions independent of the queue-lengths. Taking the $\mu_n^\prime[k]$ decisions to be the service decisions corresponding to this stationary randomized policy, we have
\begin{align*}
\mathbb{E} \{ \lambda_r(1-\alpha_r) - \mu_r^\prime[k] | \myp{k} \} &< - \epsilon, \\
\mathbb{E} \{ \beta_g - \mu_g^\prime[k] | \myp{k} \} &< - \epsilon, \\
\mathbb{E} \{ \lambda_b - \mu_b^\prime[k] | \myp{k} \} &< - \epsilon,
\end{align*} 
where $\epsilon > 0$. Thus, the previous drift inequality can the simplified as
\begin{equation*}
\Delta(\myp{k}) \leqslant B - 2 \epsilon \sum_{n \in \mathcal{N}} \Phi_n[k],
\end{equation*}
where we have used definition \eqref{def}. Now, as before, the queue evolves as a Markov chain. The drift inequality above gives negative drift but for a finite set of queue-lengths. Hence, using Theorem 2 from \cite{definition}, the network is stable.

\section*{Appendix F}

\textbf{Proof of Lemma 1:}

We are given the Lyapunov function $L(\boldsymbol{Q}(t)) = \sum_{i=1}^{N} Q_{i}^{2}(t)$. Squaring both sides of equation (\ref{Q_evolution}), and using the fact that $(max\{a,0\})^{2} \leqslant a^{2}$, we have
\[
Q_{i}^{2}(t+1) - Q_{i}^{2}(t) \leqslant \mu_{i}^{2}(t) + A_{i}^{2}(t) + 2 Q_{i}(t)(A_{i}(t)-\mu_{i}(t)).
\]
Summing over all $i$ and taking $\mathbb{E}\{\cdot\ |\ \boldsymbol{Q}(t)\}$,
\begin{align}
\Delta(\boldsymbol{Q}(t)) & \leqslant \mathbb{E}\{U(t)\ |\ \boldsymbol{Q}(t)\} + 2 \sum_{i=1}^{N} Q_{i}(t)(\lambda_{i} - \mathbb{E}\{\mu_{i}(t)\ |\ \boldsymbol{Q}(t)\})
\end{align}
where $U(t)$ is as defined in \eqref{U_t_def}. Since, $\mathbb{E}\{A_{i}^{2}(t)\} < \infty$ and $\mathbb{E}\{X_{ij}^{2}(t)\} < \infty$ for each $i$ and $j$,
\[
B \triangleq \sum_{i=1}^{N}\ \mathbb{E}\{U(t)\ |\ \boldsymbol{Q}(t)\}
\]
is bounded and independent of $\boldsymbol{Q}(t)$. Hence, we are done.



%

\end{document}